\documentclass[11pt]{article}
\usepackage{mathrsfs}
\usepackage{latexsym,amsmath,amssymb,amsfonts,epsfig}
\usepackage{color,colordvi,amscd,indentfirst,graphics}
\allowdisplaybreaks

\newtheorem{theorem}{Theorem}[section]

\newtheorem{claim}[theorem]{Claim}

\def\qed{\hfill \rule{4pt}{7pt}}
\def\pf{\noindent {\it Proof.} }

\begin{document}

\title{A note on spanning trees of connected $K_{1,t}$-free graphs whose stems have a few leaves}
\author{Pham Hoang Ha\footnote{E-mail address: ha.ph@hnue.edu.vn.(Corresponding author)}\\
Department of Mathematics\\
Hanoi National University of Education\\
136 XuanThuy Street, Hanoi, Vietnam\\
\medskip\\
Dang Dinh Hanh\footnote{E-mail address: ddhanhdhsphn@gmail.com.}\\
Department of Mathematics\\
Hanoi Architectural University\\
Km10 NguyenTrai Street, Hanoi, Vietnam\\}
\date{}

\maketitle{}

\bigskip

\begin{abstract}
Let $T$ be a tree, a vertex of degree one is called a leaf. The set of leaves of $T$
is denoted by $Leaf(T)$. The subtree $T-Leaf(T)$ of $T$ is called the stem of
$T$ and denoted by $Stem(T).$ In this note, we give a sharp sufficient condition to show that a $K_{1,t}-$free graph has a spanning tree whose stem has a few leaves. By applying the main result, we give improvements of previous related results.
\end{abstract}




\section{Introduction}

In this note, we only consider finite simple graphs. Let $G$ be a
graph with vertex set $V(G)$ and edge set $E(G)$. For any vertex
$v\in V(G)$, we use $N_G(v)$ and $\deg_G(v)$ to denote the set of neighbors of $v$ and the
degree of $v$ in $G$, respectively. For any $X\subseteq V(G)$, we denote by
$|X|$ the cardinality of $X$. We define $G-uv$ to be the graph obtained from $G$ by
deleting the edge $uv\in E(G)$, and $G+uv$ to be the graph obtained
from $G$ by adding an edge $uv$ between two non-adjacent vertices
$u$ and $v$ of $G$. For two vertices $u$ and $v$ of $G$, the distance between $u$ and $v$ in $G$ is denoted 
by $d_{G}(u, v)$.\\
For an integer $m\geqslant 2,$ let $\alpha^{m}(G)$ denote the number defined by
$$\alpha^{m}(G)=\max\{ |S|:S\subset V(G),d_{G}(x,y)\geqslant m\,\quad\text{for all vertices }\,x,y\in S\}.$$
For an integer $p\geqslant 2$, we define
$$\sigma_p^{m}(G)=\min\left\lbrace \sum_{a\in S}\deg_{G}(a):S\subset V(G),|S|=p, d_{G}(x,y)\geqslant m\quad\text{for}\;x,y\in S\right\rbrace.$$
For convenience, we define $\sigma^{m}_{p}(G)=+\infty$ if $\alpha^{m}(G)<p$. We note that,  $\alpha^{2}(G)$ is often written $\alpha(G)$, which  is the independent number of $G,$ and $\sigma_p^{2}(G)$ is often written $\sigma_{p}(G)$, which is the minimum degree sum of $p$ independent vertices.  

Let $T$ be a tree, a vertex of degree one is called a leaf. The set of leaves of $T$
is denoted by $Leaf(T)$. The subtree $T-Leaf(T)$ of $T$ is called the $stem$ of $T$ and is denoted by $Stem(T)$. A tree having at most $l$ leaves is called $l-$ended tree and a stem having at most $l$ leaves is  called $l-$ended stem. 
There are several well-known conditions (such as the independence
number conditions and the degree sum conditions) ensuring that a
graph $G$ contains a spanning tree with a bounded number of leaves
or branch vertices (see the survey paper~\cite{OY} and the
references cited therein for details). Win~\cite{Wi} obtained a
sufficient condition related to the independence number for
$k$-connected graphs, which confirms a conjecture of Las
Vergnas~\cite{LV}. Broersma and Tuinstra~\cite{BT} gave a degree
sum condition for a connected graph to contain a spanning tree with
at most $l$ leaves.
\begin{theorem}\label{theo1.1}{\rm (Win~\cite{Wi})}
	Let $G$ be a $k$-connected graph and let $l\geq 2$.  If
	$\alpha(G)\leq k+l-1$, then $G$ has a spanning tree with at most $l$
	leaves.
\end{theorem}
\begin{theorem}\label{theo1.2}{\rm (Broerma and Tuinstra~\cite{BT})}
	Let $G$ be a connected graph with $n$ vertices and let $l\geq 2$. If
	$\sigma_2(G)\geq n-l+1$, then $G$ has a spanning tree with at most
	$l$ leaves.
\end{theorem}

Recently, many researches are studied on spanning trees in connected graphs whose stems having a bounded number of leaves or branch vertices (see \cite{TZ}, \cite{KY}, \cite{KY15} and \cite{Yan} for more details). We introduce here some results on spanning trees whose stems have a few leaves. 
\begin{theorem}{\rm (Tsugaki and Zhang~\cite{TZ})}
	Let $G$ be a connected graph and let $l\geqslant 2$ be an integer. If $\sigma_{3}(G)\geqslant|G|-2l+1$, then $G$ have a spanning tree whose stem has at most $l$ leaves.
\end{theorem}
\begin{theorem}{\rm (Kano and Yan~\cite{KY})}\label{KY1}
	Let $G$ be a connected graph and let $l\geqslant 2$ be an integer. If $\sigma_{l+1}(G)\geqslant|G|-l-1$, then G has a spanning tree with $l$-ended stem.
\end{theorem}
\begin{theorem}{\rm (Kano and Yan~\cite{KY})}\label{KY2}
	Let $G$ be a connected claw-free graph and let $l\geqslant 2$ be an integer. If $\sigma_{l+1}(G)\geqslant|G|-2l-1$, then G have a spanning tree with $l$-ended stem.
\end{theorem}
On the other hand, for a positive integer $t \geq 3,$ a graph $G$ is said to be  $K_{1,t}-$ free graph if it contains no $K_{1,t}$ as an induced subgraph. If $t=3,$ the $K_{1,3}-$ free graph is also called the claw-free graph. Moreover, if the maximum degree of a graph $G$ is
denoted by $\Delta(G)$ then $G$ is nothing but a $K_{1,t}-$free graph for all $t\geq \Delta(G)+1.$ 
In this note, we would like to introduce a generalization of above theorems. We study on the spanning tree of a $K_{1,t}-$free graph whose stem has a bounded number of leaves. In particular, we prove the following result.
\begin{theorem}\label{thm-main}
	For a positive integer $t \geq 3,$ let $G$ be a connected $K_{1,t}$-free graph and let $l(\not= t-2)$ be an integer. If 
\begin{equation}\label{eq1}	
\sigma^{4}_{l+1}(G)\geqslant |G|-\lfloor \dfrac{l(t-1)}{t-2}\rfloor-1
\end{equation}
	then $G$ has a spanning tree with $l$-ended stem. Here, the notation $\lfloor r\rfloor$ stands for the biggest integer not exceed the real number $r.$
\end{theorem}
We also note that the reason why we consider $\sigma^{4}_{l+1}(G)$ is based on the following theorem of Kano and Yan.
\begin{theorem}{\rm (Kano and Yan~\cite{KY})}\label{KY}
	Let $G$ be a connected graph and let $l\geqslant 2$ be an integer. If $\alpha^{4}(G)\leq l$, then G has a spanning tree with $l$-ended stem.
\end{theorem}
By using Theorem \ref{thm-main} when $t=3,$ we have Theorem \ref{KY2}. Moreover, Theorem \ref{thm-main} is an improvement of Theorem \ref{KY1} when we consider the positive integer $t$ big enough.

 We end this section by constructing two examples to show that
the conditions of Theorem \ref{thm-main}
are sharp. Let $t, k, m$ be integers such that $t\geq 3, k\geq 2, m \geq 1$ and let $l=k(t-2)$. Let $D$ be a complete graph with $k+1$ vertices $u_1, u_2, ..., u_{k+1}$. Let $D_1, D_2, ..., D_{k(t-2)+1}$ be copies of the graph $K_m$. Let $v_1, v_2, ..., v_{k(t-2)+1}$ be vertices which are not in $D\cup D_1\cup D_2\cup \cdots \cup D_{k(t-2)+1}$. For each $i\in\{1,2,...,k\},$ join $u_i$ to all vertices of the graphs $D_{(i-1)(t-2)+1}, D_{(i-1)(t-2)+2}, ..., D_{i(t-2)} $ and join $u_{k+1}$ to all vertices of the graph $D_{k(t-2)+1}$. Join $v_j$ to all vertices of $D_j$ for all $j\in \{1;2;\cdots; k(t-2)+1\}$. Then the resulting graph $G$ is a $K_{1,t}-$free graph. Moreover, we have $|G|=k+1+(k(t-2)+1)(m+1)$ and \begin{align*}
	\sigma^4_{l+1}(G)&=\sigma^4_{k(t-2)+1}(G)=\sum\limits_{i=1}^{k(t-2)+1}\deg_G(v_i)\\ 
	&=(k(t-2)+1).m=|G|-k(t-1)-2=|G|-\lfloor \dfrac{l(t-1)}{t-2}\rfloor-2.
\end{align*}
But $G$ has no spanning tree with $l$-ended stem. Hence the condition (\ref{eq1}) is sharp.\\
On the other hand, when $l=t-2,$ let $D_1, D_2, ..., D_{l+1}$ be copies of the graph $K_m$. Let $w, v_1, v_2, ..., v_{l+1}$ be distinct vertices which are not in $ D_1\cup D_2\cup \cdots \cup D_{l+1}$. For each $i\in\{1,2,...,l+1\}$ join $v_i$ to all vertices of the graph $D_{i} $ and join $w$ to all vertices of the graphs $D_1, D_2, \cdots D_{l+1}$. The resulting graph is denoted by $H.$ Then $H$ is a $K_{1,t}-$free graph. Moreover, we may obtain that $|H|=1+(l+1)(m+1)$ and \begin{align*}
	\sigma^4_{l+1}(H)=\sum\limits_{i=1}^{l+1}\deg_{H}(v_i)=(l+1).m=|H|-l-2=|H|-\lfloor \dfrac{l(t-1)}{t-2}\rfloor-1.
\end{align*}
But $H$ has no spanning tree with $l$-ended stem. This implies that the condition $l\not=t-2$ is necessary.

\section{Proof of Theorem \ref{thm-main}}
In this section, we extend the idea of Kano and Yan in~\cite{KY} to prove
Theorem~\ref{thm-main}.

 We prove the theorem
by contradiction. Suppose to the contrary that $G$ contains no spanning tree with $l-$ended stem. We choose a maximal tree $T$ with $l-$ended stem in $G$ so that

(C) $|Leaf(T)|$ is as large as possible.\\
By the maximality of $T$, we have the following claim.
\begin{claim}\label{claim1}
	For every $v\in V(G)-V(T), N_{G}(v)\subseteq Leaf(T)\cup(V(G)-V(T))$.
\end{claim}
 Because $G$ is connected and $T$ is not a spanning tree of $G,$ there exist two vertices $v_1\in V(G)-V(T)$ and $v_2\in Leaf(T)$ such that $v_1v_2\in E(G)$.
We may obtain that $Stem(T)$ has exactly $l$ leaves. Indeed, otherwise we consider the tree $T'=T+v_1v_2.$ Then $T'$ has $l$-ended stem and $|T'|>|T|$, this implies a contradiction with the maximality of $T$. Let $\{x_1, x_2, ..., x_l\}$ be the leaf set of $Stem(T)$.
\begin{claim}\label{claim2}
	For every 
	$x_{i}(1 \leqslant i \leqslant l)$, there exists a vertex $y_{i}\in Leaf(T)$ such that $y_{i}$ is adjacent to $x_{i}$
	and $N_{G}(y_{i})\subset Leaf(T)\cup\left\{x_{i}\right\}$.
\end{claim}
\pf
	 By the maximality of $T,$ it is easy to see that for each leaf $x\in Leaf(Stem(T)),$ there exists at least a vertex $y$ in $Leaf(T)$ such that $y$ is adjacent to $x.$ Suppose that for some $1\leqslant i\leqslant l$, each leave $y_{i_j}$ of $T$ adjacent to $x_{i},$ is also adjacent to a vertex $z_{i_j}\in (Stem(T)-\left\{x_{i}\right\})$. Then we consider $T'$ to be the tree obtained from $T$ by removing the edge $y_{i_j}x_{i}$ and adding the edge $y_{i_{j}}z_{i_{j}}$. Hence $T'$ is a tree with $l-$end stem such that $|T'|=|T|$ and $Leaf(T')=Leaf(T)+\{x_i\}$, which contradicts the condition (C). Therefore, for each $x_i$, there exists a leaf $y_i\in N_G(x_i)$ such that $N_G(y_i)\cap (Stem(T)-\{x_i\})=\emptyset$. By the maximality of $T$ we also see that  $N_{G}(y_i)\cap (V(G)-V(T))=\emptyset.$ The claim holds.
\qed
\begin{claim}\label{claim3}
	For any two dictinct vertices $y,z\in \left\{v_1,y_{1},y_{2},\dots,y_{l}\right\},d_{G}\left(y,z\right)\geqslant 4$.
\end{claim}
\pf
	First, we show that $d_{G}\left(v_1,y_{i}\right)\geqslant 4$ for every $1\leqslant i\leqslant l$. Let $P_i$ be the shortest path connecting $v_1$ and $y_{i}$ in $G$.  If all the vertices of $P_i$ between $v_1$ and $y_{i}$ are contained in $Leaf(T)\cup (V(G)-V(T))\cup\left\{x_{i}\right\},$ add $P_{i}$ to $T$ (if $P_{i}$ passes through $x_{i}$, we just add the segment of $P_{i}$ between $v_1$ and $x_{i})$ and remove the edges of $T$ joining $V(P_{i}\cap Leaf(T))$ to $V(Stem(T))$ except the edge $y_{i}x_{i}.$ The resulting tree is denoted by $T'.$ Then $T'$ is a tree in $G$ with $l$-ended stem  and $|T'|>|T|$, which contradicts to the maximality of $T.$ Hence we may choose the vertex $s$ in $Stem(T)\cap P_i$ such that it is nearest to $v_1$ in $P_i$. If $s=x_j$ for some $1\leq j\leq l$, then we add the segment of $P_{i}$ between $v_1$ and $x_j$ (which is denoted by $Q$) to $T$ and remove the edges of $T$ joining $V(Q)\cap Leaf(T)$ to $V(Stem(T))$ except $x_jy_j$. Hence the resulting tree has $l$-ended stem and its order is greater then $|T|$, contradicting the maximality of $T$. Thus $s\in Stem(T)-\{x_1,...,x_l\}$. By Claims \ref{claim1} and \ref{claim2}, we have $d_G(v_1,s)\geq 2$, $d_G(s,y_i)\geq 2.$ Therefore we conclude that $d_G(v_1,y_i)=d_G(v_1,s)+d_G(s,y_i)\geq 4$. 
	
	\indent Next, we show that $d_{G}(y_{i},y_{j})\geqslant 4$ for all $1\leqslant i<j\leqslant l$. Let $P_{ij}$ be the shortest path connecting $y_{i}$ and $y_{j}$ in $G$. We note that if $P_{ij}$ passes through $x_{i}$ (or $x_{j}$), then $y_{i}x_{i}\in E(P_{ij})$ (or $y_{j}x_{j}\in E(P_{ij})$), respectively. We consider following two cases.\\
	{\it Case 1.} All vertices of $P_{ij}$ between $y_{i}$ and $y_{j}$ are contained in $Leaf(T)\cup(V(G)-V(T))\cup\left\{x_{i},x_{j}\right\}$. Then we add $P_{ij}$ to $T$ and remove the edges of $T$ joining $V(P_{ij }\cap Leaf(T))$ to $V(Stem(T))$ except the edges $y_{i}x_{i}$ and $y_{j}x_{j}$. Hence the resulting graph has exactly a cycle, which contains an edge $e$ of $Stem(T)$ incident with a branch vertex in $stem(T)$. By removing the edge $e$ and by adding an edge $v_1v_2$, we have a resulting tree $T'$ with $l$-ended stem of order greater than $|T|$, which contradicts the maximality of $T$. \\	
	{\it Case 2.} There exists a vertex $s \in P_{ij}\cap (Stem(T)-\{x_i,x_j\}).$ Then $d_G(y_i,s)\geq 2, d_G(s,y_j)\geq 2$ by Claim \ref{claim2}. This concludes that $d_G(y_i,y_j)=d_G(y_i,s)+d_G(s,y_j)\geq 4$. \\
	So the assertion of the claim holds.
\qed

Denote $Y=\left\{y_{1},y_2,\dots,y_{l}\right\}$. Since Claims \ref{claim1}-\ref{claim3}, we have
\begin{center}
	$N_{G}(v_1)\subseteq (V(G)-V(T)-\left\{v_1\right\})\cup (N_{G}(v_1)\cap(Leaf(T)-Y)),$\\
	$\displaystyle\bigcup_{i=1}^{l}N_{G}(y_{i})\subseteq (Leaf(T)-Y-N_{G}(v_1))\cup\left\{x_{1},\dots,x_{l}\right\}.$
\end{center}
Hence by setting $q=|N_{G}(v_1)\cap(Leaf(T)-Y)|$, we obtain
\begin{equation*}
	\begin{split}
		\deg_{G}(v_1)+\displaystyle\sum\limits_{i=1}^{l}\deg_{G}(y_{i})&\leqslant (|G|-|T|-1+q)+(|Leaf(T)|-l-q)+l\\
		&=|G|-|Stem(T)|-1.
	\end{split}
\end{equation*}
On the other hand, by the assumption of Theorem \ref{thm-main} we have $$|G|-\lfloor \dfrac{l(t-1)}{t-2}\rfloor-1\leq \sigma^4_{l+1}(G)\leq \deg_{G}(v_1)+\displaystyle\sum\limits_{i=1}^{l}\deg_{G}(y_{i}).$$
Therefore we obtain $|Stem(T)|\leq \lfloor \dfrac{l(t-1)}{t-2}\rfloor$. By combining with  $|Leaf(Stem(T))|=l$ we conclude that 
\begin{align} \label{1}
	|Stem(Stem(T))|\leq \lfloor \dfrac{l}{t-2}\rfloor.
\end{align}
\begin{claim}\label{claim5}
	$N_G(v_2)\cap \{x_1,x_2,...,x_l\}=\emptyset$.
\end{claim}
\pf
	Suppose the assertion of the claim is false. Then there exists a vertex $z\in N_G(v_2)\cap \{x_1,x_2,...,x_l\}$. Remove the edge of $T$ joining $v_2$ to $Stem(T)$ except $zv_2$ and add the edges $zv_2, v_1v_2.$ Hence the resulting tree has $l$-ended stem and $|T'|>|T|$, this contradicts to the condition (T1). Claim \ref{claim5} is proved.
\qed

Now, we complete the proof of Theorem \ref{thm-main} by considering following two steps.

{\bf Step 1}. $|Stem(Stem(T))|=1$. 

We assume that $Stem(Stem(T))=\{u\}$. It follows from $t\geq 3$ and (\ref{1}) that $l\geq t-2.$ But $l \not= t-2$ we have $l \geq t- 1.$ By combining with Claims \ref{claim3}- \ref{claim5}, $G$ induced a $K_{1,t}$ subgraph with the vertex set $\{u, x_1, x_2,...,x_{t-1},v_2\}$, this gives a contradiction.

{\bf Step 2}. $|Stem(Stem(T))|\geq 2$. 

By Claim \ref{claim5}, there exists a vertex $v_3\in N_G(v_2)\cap Stem(Stem(T))$.\\
Now, we conclude that $|N_T(v_3)\cap\{x_1,x_2,...,x_l\}|<t-2$. Indeed, otherwise, without loss of generality we may assume $x_1,x_2,...,x_{t-2}\in N_G(v_3)$. Setting $s\in Stem(Stem(T))\cap N_T(v_3)$. We consider the subgraph with the vertex set $\{v_3, v_2, s, x_1, x_2,...,x_{t-2}\}$ in $G.$ By combining with the fact that $G$ is $K_{1,t}-$free we have following two cases.\\
{\it Case 1.}  $sv_2\in E(G)$. This implies that the tree  $T'=T+sv_2+v_2v_1-sv_3$ has $l$-ended stem and $|T'|>|T|$, this contradicts to the maximality of $G$.\\
{\it Case 2.} $x_js \in E(G)$ for some $j \in \{1; ...; t-2\}.$ Then we consider the tree  $T'=T+x_js+x_jv_3+v_2v_1-sv_3.$ Hence $T'$ has $l$-ended stem and $|T'|>|T|$, this also contradicts to the maximality of $G$.\\
Therefore $|N_T(v_3)\cap\{x_1;x_2;...;x_l\}|<t-2$. \\
On the other hand, since $|Leaf(Stem(T))|=l$ and $|Stem(Stem(T))|\leq \lfloor \dfrac{l}{t-2}\rfloor$ there exists a vertex $u\in Stem(Stem(T))$ such that $|N_T(u)\cap\{x_1,x_2,...,x_l\}|\geq t-1$. Without loss of generality, we may assume $x_1,x_2,...,x_{t-1}\in N_T(u)$. Set $s\in Stem(Stem(T))\cap N_T(u)$. By repeating the same arguments as in {\it Case 2} we may conclude that $x_js \not\in E(G)$ for all $j \in \{1; ...; t-1\}.$ Then, $G$ induces a $K_{1,t}$ subgraph with vertex set $\{u, s, x_1,x_2,...,x_{t-1}\}$. This gives a contradiction with the assumption of Theorem \ref{thm-main}.\\
Therefore we complete the proof of Theorem \ref{thm-main}.

\bigskip

{\bf Acknowledgements.} The research is supported by the NAFOSTED Grant of Vietnam (No. 101.04-2018.03).


\begin{thebibliography}{99}
	
	\bibitem{BT}
	H. Broersma and H. Tuinstra, \emph{Independence trees and Hamilton cycles,}
	{ J. Graph Theory} {\bf 29} (1998) 227--237.

	
	\bibitem{KY} 
	\textsc{M . Kano and Z. Yan }, 
	\emph{Spanning trees whose stems have at most $k$ leaves }, Ars Combin., \textbf{CXIVII} (2014), 417-424.
	
	\bibitem{KY15} 
	\textsc{M . Kano and Z. Yan }, 
	\emph{Spanning trees whose stems are spiders}, Graphs Combin., \textbf{31} (2015), no. 6, 1883-1887.
	
	\bibitem{LV}
	M. Las Vergnas, \emph{Sur une propriet\'{e} des arbres maximaux dans un
	graphe,} { C. R. Acad. Sci. Paris Ser. A} {\bf 272} (1971)
	1297--1300.
	
	\bibitem{OY}
	\textsc{K. Ozeki and T. Yamashita}, 
	\emph{Spanning trees: A survey}, Graphs Combin., \textbf{22} (2011), 1-26. 
	
	\bibitem{TZ}
	\textsc{M. Tsugaki and Y. Zhang}, \emph{Spanning trees whose stems have a few leaves}, Ars Combin., \textbf{CXIV} (2014), 245-256.
	
	\bibitem{Yan}
	\textsc{Z. Yan}, \emph{Spanning trees whose stems have a bounded number of branch vertices,}
	Discussiones Mathematicae Graph Theory, \textbf{36} (2016), 773-778.
	
	\bibitem{Wi}
	S. Win, \emph{On a conjecture of Las Vergnas concerning certain spanning
	trees in graphs,} { Results Math.} {\bf 2} (1979) 215--224.
	
\end{thebibliography}
\end{document}